\documentclass[12pt,twoside]{article}
\usepackage[english]{babel}
\usepackage[latin1]{inputenc}
\usepackage{amsmath}
\usepackage{amssymb,amsfonts}
\usepackage{graphicx}
\usepackage{times,amssymb,amscd}
\newtheorem{thm}{Theorem}[section]

\newtheorem{lem}[thm]{Lemma}
\newtheorem{prop}[thm]{Proposition}
\newtheorem{defn}[thm]{Definition}
\newtheorem{rem}[thm]{Remark}
\numberwithin{equation}{section}

\newcommand{\bA}{\mathbf{A}}

\newcommand{\bE}{\mathbf{E}}

\newcommand{\bH}{\mathbf{H}}

\newcommand{\bL}{\mathbf{L}}

\newcommand{\bR}{\mathbf{R}}
\newcommand{\bS}{\mathbf{S}}
\newcommand{\bV}{\mathbf{V}}

\newcommand{\be}{\mathbf{e}}

\newcommand{\bT}{\mathbf{T}}

\newcommand{\bt}{\mathbf{t}}

\newcommand{\cS}{\mathcal{S}}

\newcommand{\EUC}{\mathbf E^3}
\newcommand{\SPH}{\bS^3}
\newcommand{\HYP}{\bH^3}
\newcommand{\SXR}{\bS^2\!\times\!\bR}
\newcommand{\HXR}{\bH^2\!\times\!\bR}
\newcommand{\SLR}{\widetilde{\bS\bL_2\bR}}
\newcommand{\NIL}{\mathbf{Nil}}
\newcommand{\SOL}{\mathbf{Sol}}

\begin{document}
\pagestyle{myheadings}
\markboth{\centerline{Jen\H o Szirmai}}
{Bisector surfaces and circumscribed spheres of tetrahedra $\dots$}
\title
{Bisector surfaces and circumscribed spheres of tetrahedra derived by translation curves in $\SOL$ geometry
}

\author{Jen\H o Szirmai \\
\normalsize Budapest University of Technology and \\
\normalsize Economics Institute of Mathematics, \\
\normalsize Department of Geometry \\
\normalsize Budapest, P. O. Box: 91, H-1521 \\
\normalsize szirmai@math.bme.hu
\date{\normalsize{\today}}}

\maketitle
\begin{abstract}
In the present paper we study the $\SOL$ geometry that is one of the eight homogeneous Thurston 3-geomet\-ri\-es.

We determine the equation of the translation-like bisector surface
of any two points. We prove, that the isosceles property of a translation triangle is not equivalent to two angles of the triangle being equal and that
the triangle inequalities do not remain valid for translation triangles in general.

Moreover, we develop a method to determine the centre and the radius of the circumscribed translation sphere of a given {\it translation tetrahedron}.

In our work we will use for computations and visualizations the projective model of $\SOL$ described by E. Moln\'ar in \cite{M97}.
\end{abstract}

\section{Introduction} \label{section1}
The Dirichlet - Voronoi (briefly $D-V$) cell is fundamental concept in geometry and crystallography. In particular, they do play important roles in the study of
ball packing and ball covering. In $3$-dimensional spaces of constant curvature the $D-V$ cells are widely investigated, but in the further Thurston geometries $\SXR$, $\HXR$,
$\NIL$, $\SOL$, $\SLR$ there are few results in this topic. Let $X$ be one of the above five geometries and $\Gamma$ is one of its discrete isometry groups.
Moreover, we distinguish two distance function types: $d^g$ is the usual geodesic distance function and $d^t$ is the translation distance function (see Section 3).
Therefore, we obtain two types of the $D-V$ cells regarding the two distance functions.

We define the Dirichlet-Voronoi cell with kernel point $K$ of a given discrete isometry group $\Gamma$:
\begin{defn}
We say that the point set
\[
\mathcal{D}(K)=\left\{ Y \in X : ~ d^i(K,Y)\le d^i(K^{\mathbf{g}},Y) ~ \text{for all} ~ \mathbf{g}\in \Gamma \right\}\subset X
\]
is the {\it Dirichlet-Voronoi cell} of $\Gamma$ around its {\it kernel point} $K$ where $d^i$ is the geodesic or translation distance function of $X$.
\end{defn}

The firs step to get the $D-V$ cell of a given point set of $X$ is the determination of the translation or geodesic-like bisector (or equidistant) surface
of two arbitrary points of $X$ because these surface types contain the faces of $D-V$ cells. 

In \cite{PSSz10}, \cite{PSSz11-1}, \cite{PSSz11-2} we studied the geodesic-like equidistant surfaces in $\SXR$, $\HXR$ and $\NIL$ geometries but there are no results concerning 
the translation-like equidistant surfaces in $\NIL$, $\SLR$ and $\SOL$ geometries. 

In the Thurston spaces can be introduced in a natural way (see \cite{M97}) translations mapping each point to any point.
Consider a unit vector at the origin. Translations, postulated at the
beginning carry this vector to any point by its tangent mapping. If a curve $t\rightarrow (x(t),y(t),z(t))$ has just the translated
vector as tangent vector in each point, then the  curve is
called a {\it translation curve}. This assumption leads to a system of first order differential equations, thus translation
curves are simpler than geodesics and differ from them in $\NIL$, $\SLR$ and $\SOL$ geometries. In $\EUC$, $\SPH$, $\HYP$, $\SXR$ and $\HXR$ geometries the translation and geodesic
curves coincide with each other.

Therefore, the translation curves also play an important role in $\NIL$, $\SLR$ and $\SOL$ geometries and often seem to be more natural in these geometries,
than their geodesic lines.  

{\it In this paper} we study the translation-like bisector surfaces of two points in $\SOL$ geometry, determine their equations and visualize them.
The translation-like bisector surfaces play an important role in the construction of the $D-V$ cells because their faces lie on bisector surfaces. The $D-V$-cells
are relevant in the study of tilings, ball packing and ball covering. E.g. if the point set is the orbit of a point - generated by
a discrete isometry group of $\SOL$ - then we obtain a monohedral $D-V$ cell decomposition (tiling) of the considered space and it is interesting to examine its
optimal ball packing and covering (see \cite{Sz13-1}, \cite{Sz13-2}).

Moreover, we prove, that the isosceles property of a translation triangle is not equivalent to two angles of the triangle being equal and that
the triangle inequalities do not remain valid for translation triangles in general.

Using the above bisector surfaces we develop a procedure to determine the centre and the radius of the circumscribed translation sphere of an arbitrary $\SOL$ tetrahedron.
This is useful to determine the least dense ball covering radius of a given periodic polyhedral $\SOL$ tiling because the tiling can be decomposed into tetrahedra.
\begin{rem}
We note here, that nowadays the $\SOL$ geometry is a widely investigated space concerning
its manifolds, tilings, geodesic and translation ball packings and probability theory
(see e.g. \cite{BT}, \cite{CaMoSpSz}, \cite{KV}, \cite{MSz}, \cite{MSz12}, \cite{MSzV}, \cite{Sz13-2} and the references given there).
\end{rem}
\section{On Sol geometry}
\label{sec:1}

In this Section we summarize the significant notions and notations of real $\SOL$ geometry (see \cite{M97}, \cite{S}).

$\SOL$ is defined as a 3-dimensional real Lie group with multiplication
\begin{equation}
     \begin{gathered}
(a,b,c)(x,y,z)=(x + a e^{-z},y + b e^z ,z + c).
     \end{gathered} \tag{2.1}
     \end{equation}
We note that the conjugacy by $(x,y,z)$ leaves invariant the plane $(a,b,c)$ with fixed $c$:
\begin{equation}
     \begin{gathered}
(x,y,z)^{-1}(a,b,c)(x,y,z)=(x(1-e^{-c})+a e^{-z},y(1-e^c)+b e^z ,c).
     \end{gathered} \tag{2.2}
     \end{equation}
Moreover, for $c=0$, the action of $(x,y,z)$ is only by its $z$-component, where $(x,y,z)^{-1}=(-x e^{z}, -y e^{-z} ,-z)$. Thus the $(a,b,0)$ plane is distinguished as a {\it base plane} in
$\SOL$, or by other words, $(x,y,0)$ is normal subgroup of $\SOL$.
$\SOL$ multiplication can also be affinely (projectively) interpreted by "right translations"
on its points as the following matrix formula shows, according to (2.1):
     \begin{equation}
     \begin{gathered}
     (1,a,b,c) \to (1,a,b,c)
     \begin{pmatrix}
         1&x&y&z \\
         0&e^{-z}&0&0 \\
         0&0&e^z&0 \\
         0&0&0&1 \\
       \end{pmatrix}
       =(1,x + a e^{-z},y + b e^z ,z + c)
       \end{gathered} \tag{2.3}
     \end{equation}
by row-column multiplication.
This defines "translations" $\mathbf{L}(\mathbf{R})= \{(x,y,z): x,~y,~z\in \mathbf{R} \}$
on the points of space $\SOL= \{(a,b,c):a,~b,~c \in \mathbf{R}\}$.
These translations are not commutative, in general.
Here we can consider $\mathbf{L}$ as projective collineation group with right actions in homogeneous
coordinates as usual in classical affine-projective geometry.
We will use the Cartesian homogeneous coordinate simplex $E_0(\be_0)$, $E_1^{\infty}(\be_1)$, \ $E_2^{\infty}(\be_2)$, \ 
$E_3^{\infty}(\be_3), \ (\{\be_i\}\subset \bV^4$ \ $\text{with the unit point}$ $E(\be = \be_0 + \be_1 + \be_2 + \be_3 ))$
which is distinguished by an origin $E_0$ and by the ideal points of coordinate axes, respectively.
Thus {$\SOL$} can be visualized in the affine 3-space $\bA^3$
(so in Euclidean space $\bE^3$) as well.

In this affine-projective context E. Moln\'ar has derived in \cite{M97} the usual infinitesimal arc-length square at any point
of $\SOL$, by pull back translation, as follows
\begin{equation}
   \begin{gathered}
      (ds)^2:=e^{2z}(dx)^2 +e^{-2z}(dy)^2 +(dz)^2.
       \end{gathered} \tag{2.4}
     \end{equation}
Hence we get infinitesimal Riemann metric invariant under translations, by the symmetric metric tensor field $g$ on $\SOL$ by components as usual.

It will be important for us that the full isometry group Isom$(\SOL)$ has eight components, since the stabilizer of the origin
is isomorphic to the dihedral group $\mathbf{D_4}$, generated by two involutive (involutory) transformations, preserving (2.4):
\begin{equation}
   \begin{gathered}
      (1)  \ \ y \leftrightarrow -y; \ \ (2)  \ x \leftrightarrow y; \ \ z \leftrightarrow -z; \ \ \text{i.e. first by $3\times 3$ matrices}:\\
     (1) \ \begin{pmatrix}
               1&0&0 \\
               0&-1&0 \\
               0&0&1 \\
     \end{pmatrix}; \ \ \
     (2) \ \begin{pmatrix}
               0&1&0 \\
               1&0&0 \\
               0&0&-1 \\
     \end{pmatrix}; \\
     \end{gathered} \tag{2.5}
     \end{equation}
     with its product, generating a cyclic group $\mathbf{C_4}$ of order 4
     \begin{equation}
     \begin{gathered}
     \begin{pmatrix}
                    0&1&0 \\
                    -1&0&0 \\
                    0&0&-1 \\
     \end{pmatrix};\ \
     \begin{pmatrix}
               -1&0&0 \\
               0&-1&0 \\
               0&0&1 \\
     \end{pmatrix}; \ \
     \begin{pmatrix}
               0&-1&0 \\
               1&0&0 \\
               0&0&-1 \\
     \end{pmatrix};\ \
     \mathbf{Id}=\begin{pmatrix}
               1&0&0 \\
               0&1&0 \\
               0&0&1 \\
     \end{pmatrix}.
     \end{gathered} \notag
     \end{equation}
     Or we write by collineations fixing the origin $O=(1,0,0,0)$:
\begin{equation}
(1) \ \begin{pmatrix}
         1&0&0&0 \\
         0&1&0&0 \\
         0&0&-1&0 \\
         0&0&0&1 \\
       \end{pmatrix}, \ \
(2) \ \begin{pmatrix}
         1&0&0&0 \\
         0&0&1&0 \\
         0&1&0&0 \\
         0&0&0&-1 \\
       \end{pmatrix} \ \ \text{of form (2.3)}. \tag{2.6}
\end{equation}
A general isometry of $\SOL$ to the origin $O$ is defined by a product $\gamma_O \tau_X$, first $\gamma_O$ of form (2.6) then $\tau_X$ of (2.3). To
a general point $A=(1,a,b,c)$, this will be a product $\tau_A^{-1} \gamma_O \tau_X$, mapping $A$ into $X=(1,x,y,z)$.

Conjugacy of translation $\tau$ by an above isometry $\gamma$, as $\tau^{\gamma}=\gamma^{-1}\tau\gamma$ also denotes it, will also be used by
(2.3) and (2.6) or also by coordinates with above conventions.

We remark only that the role of $x$ and $y$ can be exchanged throughout the paper, but this leads to the mirror interpretation of $\SOL$.
As formula (2.4) fixes the metric of $\SOL$, the change above is not an isometry of a fixed $\SOL$ interpretation. Other conventions are also accepted
and used in the literature.

{\it $\SOL$ is an affine metric space (affine-projective one in the sense of the unified formulation of \cite{M97}). Therefore, its linear, affine, unimodular,
etc. transformations are defined as those of the embedding affine space.}
\subsection{Translation curves}

We consider a $\SOL$ curve $(1,x(t), y(t), z(t) )$ with a given starting tangent vector at the origin $O=(1,0,0,0)$
\begin{equation}
   \begin{gathered}
      u=\dot{x}(0),\ v=\dot{y}(0), \ w=\dot{z}(0).
       \end{gathered} \tag{2.7}
     \end{equation}
For a translation curve let its tangent vector at the point $(1,x(t), y(t), z(t) )$ be defined by the matrix (2.3)
with the following equation:
\begin{equation}
     \begin{gathered}
     (0,u,v,w)
     \begin{pmatrix}
         1&x(t)&y(t)&z(t) \\
         0&e^{-z(t)}&0&0 \\
         0&0&e^{z(t)}& 0 \\
         0&0&0&1 \\
       \end{pmatrix}
       =(0,\dot{x}(t),\dot{y}(t),\dot{z}(t)).
       \end{gathered} \tag{2.8}
     \end{equation}
Thus, {\it translation curves} in $\SOL$ geometry (see \cite{MoSzi10} and \cite{MSz}) are defined by the first order differential equation system
$\dot{x}(t)=u e^{-z(t)}, \ \dot{y}(t)=v e^{z(t)},  \ \dot{z}(t)=w,$ whose solution is the following:
\begin{equation}
   \begin{gathered}
     x(t)=-\frac{u}{w} (e^{-wt}-1), \ y(t)=\frac{v}{w} (e^{wt}-1),  \ z(t)=wt, \ \mathrm{if} \ w \ne 0 \ \mathrm{and} \\
     x(t)=u t, \ y(t)=v t,  \ z(t)=z(0)=0 \ \ \mathrm{if} \ w =0.
       \end{gathered} \tag{2.9}
\end{equation}
We assume that the starting point of a translation curve is the origin, because we can transform a curve into an
arbitrary starting point by translation (2.3), moreover, unit velocity translation can be assumed :
\begin{equation}
\begin{gathered}
        x(0)=y(0)=z(0)=0; \\ \ u=\dot{x}(0)=\cos{\theta} \cos{\phi}, \ \ v=\dot{y}(0)=\cos{\theta} \sin{\phi}, \ \ w=\dot{z}(0)=\sin{\theta}; \\
        - \pi < \phi \leq \pi, \ -\frac{\pi}{2} \leq \theta \leq \frac{\pi}{2}. \tag{2.10}
\end{gathered}
\end{equation}
\begin{defn}
The translation distance $d^t(P_1,P_2)$ between the points $P_1$ and $P_2$ is defined by the arc length of the above translation curve
from $P_1$ to $P_2$.
\end{defn}
Thus we obtain the parametric equation of the the {\it translation curve segment} $t(\phi,\theta,t)$ with starting point at the origin in direction
\begin{equation}
\bt(\phi, \theta)=(\cos{\theta} \cos{\phi}, \cos{\theta} \sin{\phi}, \sin{\theta}) \tag{2.11}
\end{equation}
where $t \in [0,r] ~ r \in \bR^+$. If $\theta \ne 0$ then the system of equation is:
\begin{equation}
\begin{gathered}
        \left\{ \begin{array}{ll}
        x(\phi,\theta,t)=-\cot{\theta} \cos{\phi} (e^{-t \sin{\theta}}-1), \\
        y(\phi,\theta,t)=\cot{\theta} \sin{\phi} (e^{t \sin{\theta}}-1), \\
        z(\phi,\theta,t)=t \sin{\theta}.
        \end{array} \right. \\
        \text{If $\theta=0$ then}: ~  x(t)=t\cos{\phi} , \ y(t)=t \sin{\phi},  \ z(t)=0.
        \tag{2.12}
\end{gathered}
\end{equation}
\begin{defn} The sphere of radius $r >0$ with centre at the origin (denoted by $S^t_O(r)$) with the usual longitude and altitude parameters
$- \pi < \phi \leq \pi$,  $-\frac{\pi}{2} \leq \theta \leq \frac{\pi}{2}$, respectively, by (2.10), is specified by the equations (2.12) where $t=r$.
\end{defn}
\begin{defn}
 The body of the translation sphere of centre $O$ and of radius $r$ in the $\SOL$ space is called translation ball, denoted by $B^t_{O}(r)$,
 i.e. $Q \in B^t_{O}(r)$ iff $0 \leq d^t(O,Q) \leq r$.
\end{defn}
\begin{figure}[ht]
\centering
\includegraphics[width=13cm]{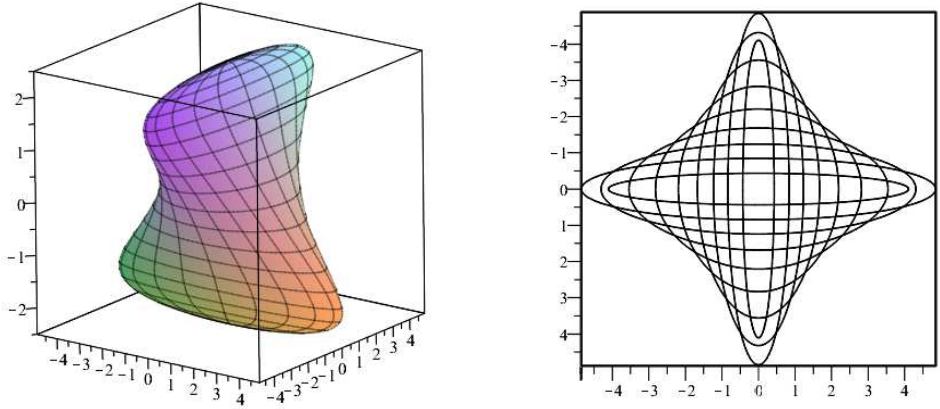}
\caption{Translation ball of radius $r=5/2$ and its plane sections parallel to [x,y] coordinate plane in $\SOL$ space}
\label{}
\end{figure}
In \cite{Sz13-2} we proved the volume formula of the translation ball $B^t_{O}(r)$ of radius $r$:
\begin{thm}
\begin{equation}
\begin{gathered}
Vol(B^t_{O}(r))=\int_{V} \mathrm{d}x ~ \mathrm{d}y ~ \mathrm{d}z = \\ = \int_{0}^{r} \int_{-\frac{\pi}{2}}^{\frac{\pi}{2}} \int_{-\pi}^{\pi} \frac{\cos{\theta}}{\sin^2{\theta}}
(e^{\rho \sin{\theta}}+e^{-\rho \sin{\theta}}-2) \ \mathrm{d}\phi \ \mathrm{d}\theta \ \mathrm{d}\rho =  \\
= 4 \pi \int_{0}^{r} \int_{-\frac{\pi}{2}}^{\frac{\pi}{2}} \frac{\cos{\theta}}{\sin^2{\theta}}
(\cosh(\rho \sin{\theta})-1) \ \mathrm{d}\theta \ \mathrm{d}\rho. \notag
\end{gathered}
\end{equation}
\end{thm}
An easy power series expansion with substitution $\rho \sin{\theta}=:z$
can also be applied, no more detailed.
From the equation of the translation spheres $S^t_O(r)$ (see (2.12)) it follows that the plane sections
of following spheres, given by parameters $\theta$ and $r$, parallel to $[x,y]$ plane are ellipses by the equations
(see Fig.~1, $r=5/2$):
\begin{equation}
\begin{gathered}
\frac{x^2}{k_1^2}+\frac{y^2}{k_2^2}=1 \ \mathrm{where} \\
k_1^2=(-\cot{\theta} (e^{-r \sin{\theta}}-1))^2, \ \ \ k_2^2=(\cot{\theta} (e^{r \sin{\theta}}-1))^2. \tag{2.13}
\end{gathered}
\end{equation}
\section{Translation-like bisector surfaces}
One of our further goals is to examine and visualize the Dirichlet-Voronoi cells of $\SOL$ geometry. In order to get $D-V$ cells we have to determine its "faces" that are
parts of bisector (or equidistant) surfaces of given point pairs. The definition below comes naturally:
\begin{defn}
The equidistant surface $\cS_{P_1P_2}$ of two arbitrary points $P_1,P_2 \in \SOL$ consists of all points $P'\in \SOL$,
for which $d^t(P_1,P')=d^t(P',P_2)$.
\end{defn}
\begin{figure}[ht]
\centering
\includegraphics[width=12cm]{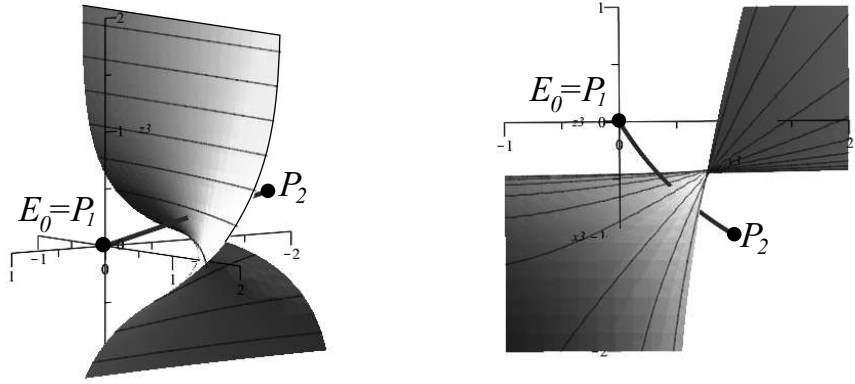}
\caption{Translation-like bisector (equidistant surface) with $P_1=(1,0,0,0)$ and $P_2=(1,-1,1,1/2)$.}
\label{pic:surf1}
\end{figure}
It can be assumed by the homogeneity of $\SOL$ that the starting point of a given translation curve segment is $E_0=P_1=(1,0,0,0)$.
The other endpoint will be given by its homogeneous coordinates $P_2=(1,a,b,c)$. We consider the translation curve segment $t_{P_1P_2}$ and determine its
parameters $(\phi,\theta,t)$ expressed by the real coordinates $a$, $b$, $c$ of $P_2$. We obtain directly by equation system (2.12) the following Lemma (see \cite{Sz17}):
\begin{lem}
\begin{enumerate}
\item Let $(1,a,b,c)$ $(b,c \in \bR \setminus \{0\}, a \in \bR)$ be the homogeneous coordinates of the point $P \in \SOL$. The paramerters of the
corresponding translation curve $t_{E_0P}$ are the following
\begin{equation}
\begin{gathered}
\phi=\mathrm{arccot}\Big(-\frac{a}{b} \frac{\mathrm{e}^{c}-1}{\mathrm{e}^{-c}-1}\Big),~\theta=\mathrm{arccot}\Big( \frac{b}{\sin\phi(\mathrm{e}^{c}-1)}\Big),\\
t=\frac{c}{\sin\theta}, ~ \text{where} ~ -\pi < \phi \le \pi, ~ -\pi/2\le \theta \le \pi/2, ~ t\in \bR^+.
\end{gathered} \tag{3.1}
\end{equation}
\item Let $(1,a,0,c)$ $(a,c \in \bR \setminus \{0\})$ be the homogeneous coordinates of the point $P \in \SOL$. The parameters of the
corresponding translation curve $t_{E_0P}$ are the following
\begin{equation}
\begin{gathered}
\phi=0~\text{or}~  \pi, ~\theta=\mathrm{arccot}\Big( \mp \frac{a}{(\mathrm{e}^{-c}-1)}\Big),\\
t=\frac{c}{\sin\theta}, ~ \text{where}  ~ -\pi/2\le \theta \le \pi/2, ~ t\in \bR^+.
\end{gathered} \tag{3.2}
\end{equation}
\item Let $(1,a,b,0)$ $(a,b \in \bR)$ be the homogeneous coordinates of the point $P \in \SOL$. The paramerters of the 
corresponding translation curve $t_{E_0P}$ are the following
\begin{equation}
\begin{gathered}
\phi=\arccos\Big(\frac{x}{\sqrt{a^2+b^2}}\Big),~  \theta=0,\\
t=\sqrt{a^2+b^2}, ~ \text{where}  ~ -\pi < \phi \le \pi, ~ t\in \bR^+.~ ~ \square
\end{gathered} \tag{3.3}
\end{equation}
\end{enumerate}
\end{lem}
{\it In order to determine the translation-like bisector surface $\cS_{P_1P_2}(x,y,z)$ of two given point $E_0=P_1=(1,0,0,0)$ and $P_2=(1,a,b,c)$
we define {translation} $\bT_{P_2}$ as elements of the isometry group of $\SOL$, that
maps the origin $E_0$ onto $P$} (see Fig.~2), moreover let $P_3=(1,x,y,z)$ a point in $\SOL$ space.

This isometrie $\bT_{P_2}$ and its inverse (up to a positive determinant factor) can be given by:
\begin{equation}
\bT_{P_2}=
\begin{pmatrix}
1 & a & b & c \\
0 & \mathrm{e}^{-c} & 0 & 0 \\
0 & 0 & \mathrm{e}^{c} & 0 \\
0 & 0 & 0 & 1
\end{pmatrix} , ~ ~ ~
\bT_{P_2}^{-1}=
\begin{pmatrix}
1 & -a\mathrm{e}^{c} & -b\mathrm{e}^{-c} & -c \\
0 & \mathrm{e}^{c} & 0 & 0 \\
0 & 0 & \mathrm{e}^{-c} & 0 \\
0 & 0 & 0 & 1
\end{pmatrix} , \tag{3.4}
\end{equation}
and the images $\bT^{-1}_{P_2}(P_i)$ of points $P_i$ $(i \in \{1,2,3\})$ are the following (see also Fig.~2):
\begin{equation}
\begin{gathered}
\bT^{-1}_{P_2}(P_1=E_0)=P_1^2=(1,-x\mathrm{e}^{z},-y\mathrm{e}^{-z},-z),~ \bT^{-1}_{P_2}(P_2)=E_0=(1,0,0,0), \\
\bT^{-1}_{P_2}(P_3)=P_3^2=(1,(x-a)\mathrm{e}^{c},(y-b)\mathrm{e}^{-c},(z-c). \tag{3.5}
\end{gathered}
\end{equation}
It is clear that $P_3=(1,x,y,z) \in \cS_{P_1P_2} ~ \text{iff} ~ d^t(P_1,P_3)=d^t(P_3,P_2) \Rightarrow d^t(P_1,P_3)=d^t(E_0,P_3^2)$ where
$P_3^2=\bT^{-1}_{P_2}(P_3)$ (see (3.4), (3.5)).

This method leads to
\begin{lem}
The implicit equation of the equidistant surface $\cS_{P_1P_2}(x,y,z)$ of two points $P_1=(1,0,0,0)$ and $P_2=(1,a,b,c)$ in $\SOL$ space (see Fig.~2,3):
\begin{enumerate}
\item $c \ne 0$
\begin{equation}\label{bis1}
\begin{gathered}
z\ne 0, c~:~
\frac{|c-z|}{|\mathrm{e}^{c}-\mathrm{e}^z|}\sqrt{(a-x)^2 \mathrm{e}^{2(c+z)}+(\mathrm{e}^{c}-\mathrm{e}^{z})^2+(b-y)^2}=\\
=\frac{|z|}{|\mathrm{e}^{z}-1|}\sqrt{x^2 \mathrm{e}^{2z}+(\mathrm{e}^{z}-1)^2+y^2},\\
z=c~:~\sqrt{(x-a)^2\mathrm{e}^{2c}+(y-b)^2\mathrm{e}^{-2c}}
=\frac{|z|}{|\mathrm{e}^{z}-1|}\sqrt{x^2 \mathrm{e}^{2z}+(\mathrm{e}^{z}-1)^2+y^2},\\
z=0~:~ \frac{|c|}{|\mathrm{e}^{c}-1|}\sqrt{(a-x)^2 \mathrm{e}^{2c}+(\mathrm{e}^{c}-1)^2+(b-y)^2}=\sqrt{x^2+y^2},
\end{gathered} \tag{3.6}
\end{equation}
\item $c=0$
\begin{equation}\label{bis1}
\begin{gathered}
z\ne 0~:~
\frac{|z|}{|\mathrm{e}^z-1|}\sqrt{(a-x)^2 \mathrm{e}^{2z}+(\mathrm{e}^{z}-1)^2+(b-y)^2}=\\
=\frac{|z|}{|\mathrm{e}^{z}-1|}\sqrt{x^2 \mathrm{e}^{2z}+(\mathrm{e}^{z}-1)^2+y^2}\Leftrightarrow \mathrm{e}^{2z}a(a-2x)+b(b-2y)=0,\\
z=0~:~ \sqrt{(x-a)^2+(y-b)^2}=\sqrt{x^2+y^2}\Leftrightarrow xa+yb-\frac{a^2+b^2}{2}.~\square
\end{gathered} \tag{3.7}
\end{equation}
\end{enumerate}
\end{lem}
\begin{figure}[ht]
\centering
\includegraphics[width=12cm]{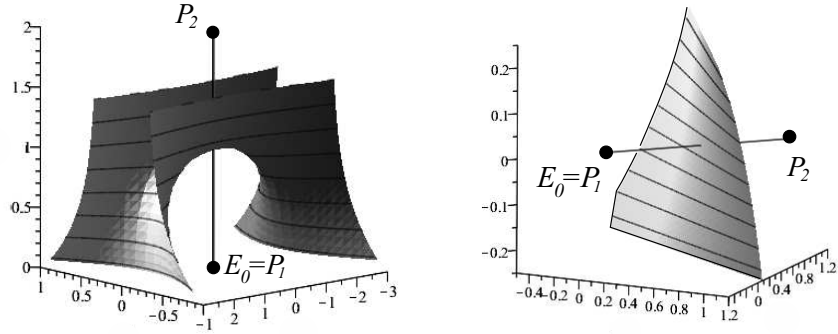}
\caption{Translation-like bisectors (equidistant surfaces) of point pairs $(P_1,P_2)$ with coordinates $((1,0,0,0), (1,0,0,2))$ (left) and $((1,0,0,0), (1,1,1,0))$ (right) }
\label{pic:surf1}
\end{figure}
\subsection{On isosceles and equilateral translation triangles}
We consider $3$ points $A_1$, $A_2$, $A_3$ in the projective model of $\SOL$ space.
The {\it translation segments} connecting the points $A_i$ and $A_j$
$(i<j,~i,j,k \in \{1,2,3\})$ are called sides of the {\it translation triangle} $A_1A_2A_3$. The length of the side $a_k$
$(k\in \{1,2,3\})$ of a translation triangle
$A_1A_2A_3$ is the translation distance $d^t(A_i,A_j)$ between the vertices $A_i$ and $A_j$ $(i<j,~i,j,k \in \{1,2,3\}, k \ne i,j$).

Similarly to the Euclidean geometry we can define the notions of isosceles and equilateral translation triangles.

An isosceles translation triangle is a triangle with (at least) two equal sides and a triangle
with all sides equal is called an equilateral translation triangle (see Fig.~4) in the $\SOL$ space.

We note here, that if in a translation triangle $A_1A_2A_3$ e.g. $a_1=a_2$ then the bisector surface $\cS_{A_1A_2}$ contains the vertex $A_3$ (see Fig.~4).
\begin{figure}[ht]
\centering
\includegraphics[width=12cm]{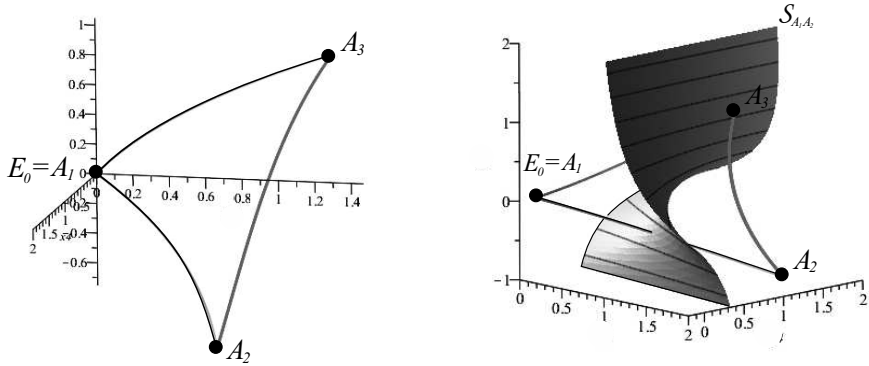}
\caption{Equilateral translation triangle with vertices $A_1=(1,0,0,0)$, $A_2=(1,2,1,-3/4)$, $A_3=(1,1,\approx 1.46717, \approx 1.04627)$ (left) and the above triangle with bisector
$\cS_{A_1A_2}$ containing the vertex $A_3$ (right).}
\label{}
\end{figure}

In the Euclidean space the isosceles property of a triangle is equivalent to two angles of the triangle being equal
therefore has both two equal sides and two equal angles. An equilateral triangle is a special case of an isosceles triangle having not just two, but all three sides and angles equal.
\begin{prop}
The isosceles property of a translation triangle is not equivalent to two angles of the triangle being equal in the $\SOL$ space.
\end{prop}
{\bf Proof:}~
The coordinates $y^3$, $z^3$ of the vertex
$A_3$ can be determined by the equation system $d^t(A_1,A_2)=d^t(A_1,A_3)=d^t(A_2,A_3)$, $y^3 \approx 1.46717$, $z^3 \approx 1.04627$
$(a_3=d^t(A_1,A_2)=a_2=d^t(A_1,A_3)=a_1=d^t(A_2,A_3) \approx 2.09436)$ (see Fig.~4).

The {\it interior angles} of translation triangles are denoted at the vertex $A_i$ by $\omega_i$ $(i\in\{1,2,3\})$.
We note here that the angle of two intersecting translation curves depends on the orientation of their tangent vectors. 

{\it In order to determine the interior angles of a translation triangle $A_1A_2A_3$
and its interior angle sum $\sum_{i=1}^3(\omega_i)$,}
we apply the method (we do not discuss here) developed in \cite{Sz17} using the infinitesimal arc-lenght square of $\SOL$ geometry (see (2.4)).

Our method (see \cite{Sz17}) provide the following results:
$$
\omega_1 \approx 0.94694,~\omega_2 \approx 1.04250,~ \omega_3 \approx 1.44910,~\sum_{i=1}^3(\omega_i) \approx 3.43854 > \pi.
$$
From the above results follows the statement. We note here, that if the vertices of the translation triangle lie in the $[x,y]$ plane than the
Euclidean isosceles property true in the $\SOL$ geometry, as well. ~$\square$

Using the above methods we obtain the following
\begin{lem}
The triangle inequalities do not remain valid for translation triangles in general.
\end{lem}
{\bf Proof:}~
We consider the translation triangle $A_1A_2A_3$ where $A_1=(1,0,0,0)$, $A_2=(1,-1,2,1)$, $A_3=(1,3/4,3/4,1/2)$.
We obtain directly by equation systems (3.1), (3.2), (3.3) (see Lemma 3.2 and \cite{Sz17}) the lengths of the translation segments $A_iA_j$ $(i,j \in \{1,2,3\}$, $i<j)$:
\begin{equation}
\begin{gathered}
d^t(A_1A_2) \approx 2.20396,~d^t(A_1A_3) \approx 1.22167, ~ d^t(A_2A_3) \approx 3.74623, \\ \text{therefore}
~ d^t(A_1A_2)+d^t(A_1A_3) < d^t(A_2A_3). ~\square \tag{3.8}
\end{gathered}
\end{equation}
We note here that if the vertices of a translation triangle lie on the $[x,y]$ plane of the model then the
corresponding triangle inequalities are true (see (2.12) and Lemma 3.2).
\subsection{The locus of all points equidistant from three given points}
A point is said to be equidistant from a set of objects if the distances between that point and each object in the set are equal.
Here we study that case where the objects are vertices of a $\SOL$ translation triangle $A_1A_2A_3$ and determine the locus of all points that are equidistant from $A_1$, $A_2$ and $A_3$.

We consider $3$ points $A_1$, $A_2$, $A_3$ that do not all lie in the same translation curve in the projective model of $\SOL$ space.
The {\it translation segments} connecting the points $A_i$ and $A_j$ $(i<j,~i,j,k \in \{1,2,3\}, k \ne i,j$) are called sides of the {\it translation triangle} $A_1A_2A_3$.
The locus of all points that are equidistant from the vertices $A_1$, $A_2$ and $A_3$ is denoted by $\mathcal{C}$.

In the previous section we determined the equation of translation-like bisector (equidistant) surface to any two points in the $\SOL$ space.
It is clear, that all points on the locus $\mathcal{C}$ must lie on the equidistant surfaces $\cS_{A_iA_j}$, $(i<j,~i,j \in \{1,2,3\})$ therefore
$\mathcal{C}=\mathcal{S}_{A_1A_2} \cap \mathcal{S}_{A_1A_3}$ and the coordinates of each of the points of that locus and only those points must satisfy the corresponding
equations of Lemma 3.3. Thus, the non-empty point set $\mathcal{C}$ can be determined and can be visualized for any given translation triangle (see Fig.~5 and 6).
\begin{figure}[ht]
\centering
\includegraphics[width=12cm]{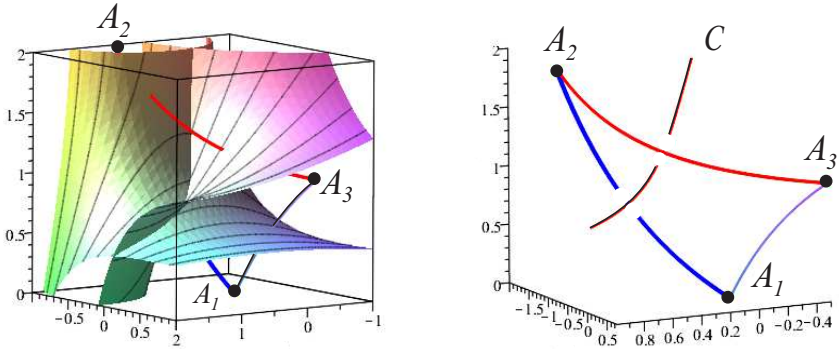}
\caption{Translation triangle with vertices $A_1=(1,0,0,0)$, $A_2=(1,2,1,-3/4)$, $A_3=(1,1,-1/2,2/3)$ with translation-like bisector surfaces $\cS_{A_1A_2}$ and $\cS_{A_1A_3}$ (left)
and a part of the locus $\mathcal{C}=\cS_{A_1A_2} \cap \cS_{A_1A_3}$ of all points equidistant from three given points $A_1$, $A_2$, $A_3$ (right).}
\label{}
\end{figure}
If the vertices of the translation triangle lie on the $[x,y]$ plane $A_1=(1,0,0,0)$, $A_2=(1,a,b,0)$, $A_3=(1,a_1,b_1,0)$ then
the parametric equation $(z\in \mathbf{R})$ of $\mathcal{C}$ is the following (see Lemma 3.3 and Fig.6):
\begin{equation}
\mathcal{C}(z)=\Big(\frac{bb_1(b-b_1)\mathrm{e}^{-2z}+a^2b_1-a_1^2b)}{2(ab_1-a_1b)},\frac{(-aa_1(a-a_1)\mathrm{e}^{2z}+ab_1^2-a_1b^2)}{2(ab_1-a_1b)},z\Big). \tag{3.9}
\end{equation}
\begin{figure}[ht]
\centering
\includegraphics[width=12cm]{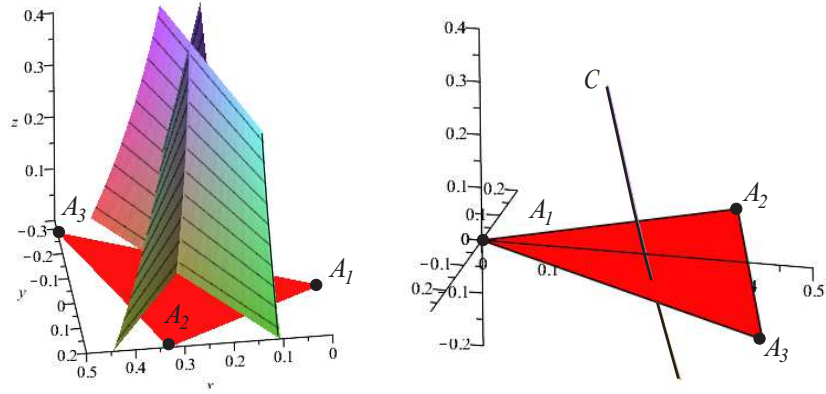}
\caption{Translation triangle with vertices $A_1=(1,0,0,0)$, $A_2=(1,1/3,1/5,$ $0)$, $A_3=(1,1/2,-2/7,0)$ with translation-like bisector surfaces $\cS_{A_1A_2}$ and $\cS_{A_1A_3}$ (left)
and a part of the locus $\mathcal{C}=\cS_{A_1A_2} \cap \cS_{A_1A_3}$ of all points equidistant from three given points $A_1$, $A_2$, $A_3$ (right).}
\label{}
\end{figure}
\subsection{Translation tetrahedra and their circumscribed spheres}
We consider $4$ points $A_1$, $A_2$, $A_3$, $A_4$ in the projective model of $\SOL$ space (see Section 2).
These points are the vertices of a {\it translation tetrahedron} in the $\SOL$ space if any two {\it translation segments} connecting the points $A_i$ and $A_j$
$(i<j,~i,j \in \{1,2,3,4\}$) do not have common inner points and any three vertices do not lie in a same translation curve.
Now, the translation segments $A_iA_j$ are called edges of the translation tetrahedron $A_1A_2A_3A_4$.
\begin{figure}[ht]
\centering
\includegraphics[width=12cm]{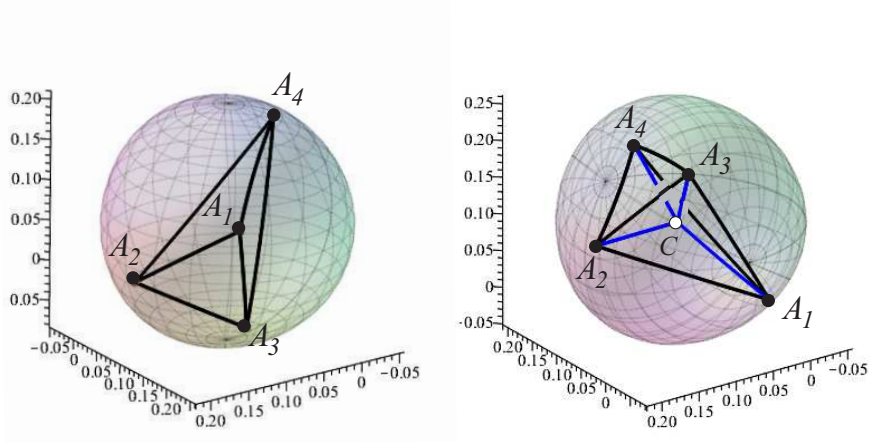}
\caption{Translation tetrahedron with vertices $A_1=(1,0,0,0)$, $A_2=(1,\sqrt{3}/8,$ $1/8,1/40)$, $A_3=(1,1/8,\sqrt{3}/8,-1/40)$, $A_4=(1,1/20,3/20,1/5)$
and its circumscibed sphere of radius $r \approx 0.14688$ with circumcenter $C=(1,\approx 0.08198, \approx 0.10540, \approx 0.06319)$.}
\label{}
\end{figure}
The circumscribed sphere of a translation tetrahedron is a translation sphere (see Definition 2.2, (2.12) and Fig.~1) that touches each of the tetrahedron's vertices.
As in the Euclidean case the radius
of a translation sphere circumscribed around a tetrahedron $T$ is called the circumradius of $T$, and the center point of this sphere is called the circumcenter of $T$.

\begin{lem}
For any translation tetrahedron there exists uniquely a translation sphere (called the circumsphere) on which all four vertices lie.
\end{lem}
\begin{figure}[ht]
\centering
\includegraphics[width=12cm]{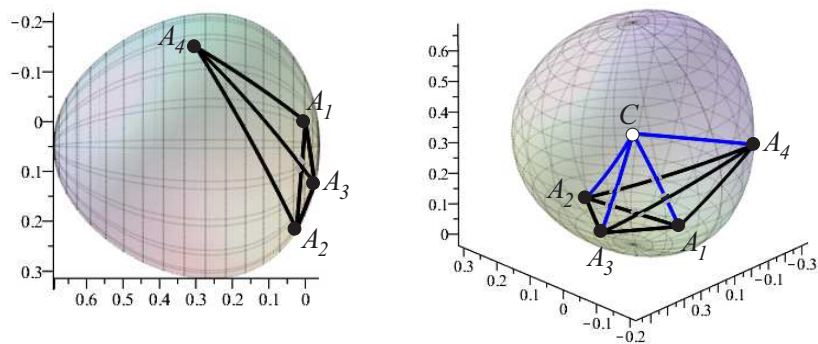}
\caption{Translation tetrahedron with vertices $A_1=(1,0,0,0)$, $A_2=(1,\sqrt{3}/8,$ $1/8,1/40)$, $A_3=(1,1/8,\sqrt{3}/8,-1/40)$, $A_4=(1,-3/20,-3/20,$ $3/10)$
and its circumscibed sphere of radius $r \approx 0.36332$ with circumcenter $C=(1,\approx 0.04904, \approx 0.17721, \approx 0.32593)$. }
\label{}
\end{figure}
{\bf Proof:}~ The Lemma follows directly from the properties of the translation distance function (see Definition 2.1 and (2.12)).
The procedure to determine the radius and the circumcenter of a given translation tetrahedron is the folowing:

The circumcenter $C=(1,x,y,z)$ of a given translation tetrahedron $A_1A_2A_3A_4$ $(A_i=(1,x^i,y^i,z^i), ~ i \in \{1,2,3,4\})$
have to hold the following system of equation:
\begin{equation}
d^t(A_1,C)=d^t(A_2,C)=d^t(A_3,C)=d^t(A_4,C), \tag{3.10}
\end{equation}
therefore it lies on the translation-like bisector surfaces $\cS_{A_i,A_j}$ $(i<j,~i,j \in \{1,2,3,4\}$) which equations are determined in Lemma 3.3.
The coordinates $x,y,z$ of the circumcenter of the circumscribed sphere around the tetrahedron $A_1A_2A_3A_4$ are obtained by the system of equation
derived from the facts:
\begin{equation}
C \in \cS_{A_1A_2}, \cS_{A_1A_3}, \cS_{A_1A_4}. \tag{3.11}
\end{equation}
Finally, we get the circumradius $r$ as the translation distance e.g. $r=d^t(A_1,C)$.

We apply the above procedure to two tetrahedra determined their centres and the
radii of their circumscribed balls that are described in  Fig.~7 and 8.~ $\square$
%

\end{document}